\documentclass[12pt]{amsart}

\usepackage{amsmath,amsfonts,amssymb,amscd}
\begin{document}
\title{Stability of extremal K\"ahler manifolds}
\author{Toshiki Mabuchi}
\address{Department of Mathematics, Graduate School of Science, Osaka University, 
Toyonaka, Osaka, 560-0043 Japan
}

\date{}
\maketitle

\thanks{\em \qquad \;\,\,Dedicated to Professor Shoshichi Kobayashi 
on his seventieth birthday}
\footnotetext{\, \, \,
To appear in Osaka Journal of Mathematics {\bf 41}(2004).}
\footnotetext{\; \; ${}^\dag{}$For this uniquness, we choose
$Z^{\Bbb C}$ (cf. Section 2) as the algebraic torus $T$.}

\section{Introduction}
In Donaldson's study  \cite{Don2} of asymptotic stability for
 polarized algebraic
manifolds $(M, L)$, 
{\it critical
metrics\/} originally defined by Zhang
\cite{Z2} (see also \cite{Luo}) are referred to as balanced
metrics and play a central role when the polarized algebraic manifolds admit K\"ahler
metrics of constant scalar curvature.    Let
$T\cong (\Bbb C^*)^k$ be an  algebraic torus in the identity component
$\operatorname{Aut}^0(M)$ of  the group of holomorphic automorphisms of
$M$. In this paper,  we define the concept of {\it critical metrics relative to\/} $T$,  
and  as an application, choosing a suitable $T$, 
we shall show that a result in \cite{M4} 
on the asymptotic approximation of critical metrics (see \cite{Don2}, \cite{Z2}) 
can be generalized 
to the case where $(M, L)$ admits
an extremal K\"ahler metric in the polarization class.
Then in our forthcoming paper \cite{M5},  we shall show that a slight modification of the
concept of stability  (see Theorem A below) allows us to obtain the asymptotic stability of extremal K\"ahler
manifolds even when the obstruction as in \cite{M4} does not vanish. 
In particular, by an argument similar to \cite{Don2}, 
an extremal K\"ahler metric in a fixed integral K\"ahler class on a projective algebraic manifold $M$ 
will be shown to be unique${}^\dag$ up to the action of the group $\operatorname{Aut}^0(M)$.

\section{Statement of results}

Throughout this paper, we fix once for all an ample holomorphic line bundle $L$
on a connected projective algebraic manifold $M$.
Let $H$ be the maximal connected linear algebraic subgroup of $\operatorname{Aut}^0(M)$,
so that $\operatorname{Aut}^0(M)/H$ is an abelian variety.
The corresponding Lie subalgebra of $H^0(M, \mathcal{O}(T^{1,0}M))$
will be denoted by $\frak h$.
For the complete linear system $|L^m|$, 
$m \gg 1$,
we consider the  Kodaira embedding 
$$
\Phi_m = \Phi_{|L_{}^{m}|} : M\; \hookrightarrow \;\Bbb P^* (V_m), \qquad m \gg
1,
$$
where $\Bbb P^* (V_m)$ denotes the set of all hyperplanes through the origin in
 $V_m := H^0(M, \mathcal{O}(L^{
m}))$. Put $N_m := \dim V_m -1$.
Let $n$ and $d$ be respectively the dimension of $M$ and the
degree of  the image $M_m := \Phi_ m
(M )$ in the projective space $\Bbb P^* (V_m)$. 
Put $W_m =\{\operatorname{Sym}^d (V_m )\}^{\otimes n+1}$.
Then to the image $M_m$ of $M$,
we can associate a nonzero element
$\hat{M}_m$ in $W^*_m$ such that 
the corresponding element $[\hat{M}_m ]$ in 
$\Bbb P^* (W_m)$ is the Chow point associated to the 
irreducible reduced algebraic cycle
$M_m$ on $\Bbb P^* (V_m)$. 
Replacing $L$ by some positive integral multiple of $L$ if necessary, we
fix an $H$-linearization of $L$, i.e.,
a lift to $L$ of the $H$-action on $M$ such that
$H$ acts on $L$ as bundle isomorphisms covering the $H$-action on $M$.
For an algebraic torus $T$ in $H$,
this naturally induces a $T$-action on $V_m$ for each $m$.
Now for each character $\chi \in \operatorname{Hom}(T, \Bbb C^*)$,
we set
$$
V(\chi )\; :=\; 
\left \{\, s  \in V_m \; ;\; t \cdot s = \chi (t ) \, s\,\text{ for all $t \in
T$}\,
\right \}.
$$
Then we have mutually distinct characters $\chi_1$, $\chi_2$, 
\dots , $\chi_{\nu_m}
\in \operatorname{Hom}(T, \Bbb C^*)$
 such that the vector space
$V_m = H^0(M, \mathcal{O}(L^m))$ is uniquely written as a direct sum
$$
V_m\; =\; \bigoplus_{k =1}^{\nu_m}\; V(\chi_k). 
\leqno{(2.1)}
$$
Put $G_m := \Pi_{k =1}^{\nu_m} \operatorname{SL}(V(\chi_k))$,
and the associated Lie subalgebra of $\operatorname{sl}(V_m)$ 
will be denoted by $\frak g_m$.
 More precisely, $G_m$ and $\frak g_m$ possibly
depend on the choice of the algebraic torus $T$,
and if necessary, we denote these by $G_m (T)$ and $\frak g_m (T)$, respectively.
The $T$-action on $V_m$ is, more precisely, a right action, while we 
regard the $G_m$-action on $V_m$
as a left action. Since $T$ is Abelian, this $T$-action on $V_m$ 
can be regarded also as a left action.

\smallskip
 The group $G_m$ acts diagonally on
$V_m$  in such a way that, for each $k$,
the $k$-th factor $\operatorname{SL}(V(\chi_k))$ of $G_m$
acts  just on the $k$-th factor 
$V(\chi_k)$ of $V_m$. This induces a natural $G_m$-action on $W_m$ and also on $W_m^*$.

\medskip
{\em Definition $2.2$.} (a) The subvariety $M_m$ of $\Bbb P^*(V_m)$ is said to be
{\it stable relative to $T$\/} or {\it semistable relative to $T$\/}, according as the orbit 
$G_m \cdot \hat{M}_m$ is closed in $W^*_m$ or the closure of $G_m \cdot \hat{M}_m$ in $W^*_m$
does not contain the origin of $W^*_m$.

\smallskip\noindent
(b) Let $\frak t_c$ denote the Lie subalgebra of the maximal compact
subgroup $T_c$ of $T$, and as a real Lie subalgebra of the complex Lie algebra $\frak t$,
we define $\frak t_{\Bbb R} := \sqrt{-1} \,\frak t_c$.

\medskip 
Take a Hermitian metric for $V_m$ such that $V (\chi_k ) \perp V(\chi_{\ell})$ if $k \neq \ell $.
Put $N_m := \dim V_m -1$ and $n_k := \dim V(\chi_k )$. We then set
$$
l (k, i) \; :=\; (i - 1)\, +\, \sum_{j=1}^{k-1} n_j,
\qquad i = 1,2,\dots, n_k; \; k = 1,2,\dots, \nu_m,
$$
where the right-hand side denotes $i-1$ in the special case  $k = 1$.
Let $\|\; \|$ denote the Hermitian norm for $V_m$ induced by the Hermitian metric.
Take a $\Bbb C$-basis $\{ s_0, s_1, \dots, s_{N_m}\}$ for $V_m$.

\medskip
{\em Definition $2.3$.}  We say that $\{s_0, s_1, \dots, s_{N_m}\}$ is an 
{\it admissible normal basis\/}  for $V_m$ if there exist positive real constants
$b_k$, $k = 1,2, \dots, \nu_m$,  and 
a $\Bbb C$-basis $\{s_{k,i}\,;\, i=1,2,\dots,n_k\}$ for $V(\chi_k )$,
 with $\Sigma_{k=1}^{\nu_m}\, n_k b_k =
N_m + 1$, such that
\begin{enumerate}
\item[(1)]  \; $s_{l (k,i)} = s_{k,i}$, \qquad $i = 1,2,\dots, n_k$;\; $k =1,2, \dots, \nu_m$;
\item[(2)]  \; $s_l \perp s_{l'}$\;  if $l \neq l'$;
\item[(3)] \; $\| s_{k, i}\|^2\;=\; b_k$, \qquad $i = 1,2,\dots, n_k$;\; $k =1,2, \dots, \nu_m$.
\end{enumerate}
Then the real vector $b := (b_1, b_2,\dots, b_{\nu_m})$ is called the {\it index\/} of
the admissible normal basis $\{s_0, s_1, \dots, s_{N_m} \}$ for $V_m$.

\medskip
We now specify a Hermitian metric on $V_m$. 
For the maximal compact
subgroup $T_c$ of $T$ above, 
let $\mathcal{S}$ be the  set $\,(\neq \emptyset )\,$ of all $T_c$-invariant K\"ahler forms in the class
$c_1(L)_{\Bbb R}$.  Let $\omega \in \mathcal{S}$,
and choose a Hermitian metric $h$ for $L$ such that $\omega = c_1(L; h)$.
Define a Hermitian metric  on $V_m$ by
$$
(s, s')_{L^2}\;\; := \int_M (s, s')^{}_{h^m} \, \omega^n, \qquad s, s'\in V_m,
\leqno{(2.4)}
$$
where $(s, s')^{}_{h^m}$ denotes the function on $M$ 
obtained as the the pointwise inner product of $s$, $s'$ by
the Hermitian metric $h^m$ on $L^m$.  
Now, let us consider the situation that $V_m$ has the  Hermitian metric (2.4).
Then
$$
V(\chi_k ) \perp  V(\chi_{\ell} ), \qquad k \neq \ell,
$$
and  define a maximal compact subgroup $(G_{m})_c$ of $G_m$ by
  $(G_{m})_c := \Pi_{k=1}^{\nu_m}\operatorname{SU}(V(\chi_k ))$. 
Again by this Hermitian metric $(\; ,\; )_{L^2}$, 
 let $\{s_0, s_1, \dots, s_{N_m}\}$ an admissible normal basis for $V_m$ 
of a given index $b$. Put
$$
E_{\omega, b } \; :=\; \sum_{i=0}^{N_m}\;  |s_i |_{h^m}^{\,2},
\leqno{(2.5)}
$$
where $|s |^{}_{h^m} := (s,s)^{}_{h^m}$ for all $s \in V_m$.
Then $E_{\omega, b}$ 
depends only on  $\omega$ and $b$. Namely, once $\omega$ and $b$ are fixed, 
$E_{\omega, b}$ is independent of the choice
of an admissible normal basis for $V(\chi_k )$ of  index $b$. 
Fix a positive integer $m$ such that $L^m$ is very ample.

\medskip
{\em Definition $2.6$.}  An element $\omega$ in $\mathcal{S}$ is called a
{\it critial metric relative to $T$\/}, if there exists an admissible normal basis 
$\{s_0, s_1, \dots, s_{N_m}\}$ for $V_m$ such that
the associated function $E_{\omega, b}$ on $M$ is constant for the index
$b$ of the admissible normal basis.
This generalizes  a {\it critical metric\/}
of Zhang \cite{Z2} (see also \cite{CGR}) who treated the case $T =\{1\}$.
If $\omega$ is a critical metric relative to $T$, then by integrating the equality (2.5) over $M$,
we see that the constant
$E_{\omega, b}$ is $(N_m +1)/c_1(L)^n[M]$.

\medskip
For the   centralizer $Z_H(T)$ of $T$ in $H$, let $Z_H(T)^0$ be its identity component.
For $m$ as above, the following generalization of a result in  \cite{Z2} is  crucial
to our study of stability:

\medskip
{\bf Theorem A.} 
{\em  The subvariety $M_m$ of $\Bbb P (V_m)$ is stable relative to $T$ 
if and only if there exists a critical metric $\omega\in \mathcal{S}$  relative to $T$.
Moreover, for a fixed index $\,b$, 
a critical metric $\omega$ in $\mathcal{S}$ relative to $T$ with  constant 
$E_{\omega, b}$   is unique up to the
action of 
$Z_H(T)^0$.}

\medskip
We now fix a maximal compact connected subgroup $K$ of $H$.
The corresponding Lie subalgebra of $\frak h$ is denoted by $\frak k$.
Let $\mathcal{S}_{K}$ denote the set of all K\"ahler forms $\omega$ in the class 
$c_1(L)_{\Bbb R}$ such that the identity component of the group of the isometries 
of $(M, \omega )$ coincides with $K$.
Then $\mathcal{S}_{K} \neq \emptyset$,
and an extremal K\"ahler metric, if any, in the class $c_1(L)_{\Bbb R}$ 
is always in  $H$-orbits of elements of $\mathcal{S}_K$.
For each $\omega\in \mathcal{S}_K$, we write 
$$
\omega \;=\; \frac{\sqrt{-1}}{2\pi}\,
\sum_{\alpha, \beta} g_{\alpha\bar{\beta}} dz^{\alpha}\wedge
dz^{\bar{\beta}}
$$
in terms of a system $(z^1, \dots, z^n)$ of holomorphic local coordinates on $M$.
let  $\mathcal{K}_{\omega}$ be the space of
all real-valued smooth functions $u$ on $M$ such that 
$\int_M u \omega^n = 0$ and that 
$$
\operatorname{grad}_{\omega}^{\Bbb
C} u\; :=\; \frac{1}{\sqrt{-1}}\, 
\sum_{\alpha, \beta} \,g^{\bar{\beta}\alpha}
\frac{\partial u}{\partial z^{\bar{\beta}}}\frac{\partial\;}{\partial z^{\alpha}}
$$
is a holomorphic vector field on $M$. Then $\mathcal{K}_{\omega}$ forms a
real Lie subalgebra of $\frak h$ by the Poisson bracket for $(M, \omega )$.
We then have the Lie algebra isomorphism
$$
\mathcal{K}_{\omega} \cong \frak k,\;\; \qquad u \,\leftrightarrow \,
\operatorname{grad}_{\omega}^{\Bbb
C} u. 
$$
For the space $C^{\infty}(M)_{\Bbb R}$ of real-valued smooth functions  on $M$, 
we consider the inner product defined by
$(u_1, u_2)_{\omega} := \int_M u_1 u_2\, \omega^n$ for
$u_1$, $u_2\in C^{\infty}(M)_{\Bbb R}$.
Let $\operatorname{pr}: 
C^{\infty}(M)_{\Bbb R} \to \mathcal{K}_{\omega}$ be the orthogonal projection.
Let $\frak z$ be the center of $\frak k$.
Then the vector field
$$
\mathcal{V} := \operatorname{grad}_{\omega}^{\Bbb
C}\operatorname{pr}(\sigma_{\omega})\in
\frak z
$$
is callled the {\it extremal K\"ahler vector field\/} of $(M, \omega )$, where
$\sigma_{\omega}$ denotes the scalar curvature of $\omega$.
Then $\mathcal{V}$ is  independent of the
choice of
$\omega$ in
$\mathcal{S}$, and satisfies $\exp (2\pi \gamma \mathcal{V}) = 1$ for some positive
integer $\gamma$ (cf. \cite{FM1}, \cite{Nakagawa}). 
Next, since we have an $H$-linearization of $L$, there exists a natural inclusion $H\subset
\operatorname{GL}(V_m)$. By passing to the Lie algebras, we obtain 
$$
\frak h\; \subset \; \frak{gl}(V_m).
$$
Take a Hermitian metric $h$ for $L$ such that the corresponding first Chern form
$c_1(L; h)$ is $\omega$.
As in \cite{M1}, (1.4.1), the infinitesimal $\frak h$-action on $L$ induces an
infinitesimal
$\frak h$-action on the  complexification $\mathcal{H}^{\Bbb C}_m$ 
of the space of all
Hermitian metrics $\mathcal{H}_m$ on the line bundle $L^m$. 
The  Futaki-Morita character
$F : \frak h \to \Bbb C$ is given by
$$
F (\mathcal{Y} )\;  := \; \frac{\sqrt{-1}}{2\pi}\int_M h^{-1}
(\mathcal{Y} h)\;
\omega^n,
$$
which is independent of the choice of $h$ (see for instance \cite{FM3}).
For the identity component $Z$ of the center of $K$, we consider
its complexification $Z^{\Bbb C}$ 
in $H$.  Then the corresponding Lie algebra is just the complexification $\frak z^{\Bbb C}$ of 
$\frak z$ above.
 We now consider the set $\Delta$ of all algebraic tori in $Z^{\Bbb C}$.  Let $T \in \Delta$.
 Put 
$$
q:= 1/m.
$$
For $\omega = c_1(L; h)\in
\mathcal{S}_K$, we consider the Hermitian metric (2.4) for $V_m$. 
We then choose an admissible normal basis  $\{s_0, s_1, \dots, s_{N_m}\}$ for $V_m$ of
index $(1,1,\dots,1)$. By the asymptotic expansion of Tian-Zelditch (cf. 
\cite{T1}, \cite{Z1}; see also \cite{Cat}) for
$m
\gg 1$, there exist real-valued smooth functions $a_k (\omega)$, $k$ = 1,2,\dots, on $M$ such that
$$
\frac{n!}{m^n}\sum_{j =0}^{N_m}\; | s_j |^{\,2}_{h^m}\;
 =\; 1 +  a_1(\omega ) q + a_2(\omega ) q^2 + \cdots\; \; .
\leqno{(2.7)}
$$
 Then $a_1(\omega ) = \sigma_{\omega}/2$ by a result of Lu \cite{Lu}. 
Let $\mathcal{Y}\in\frak t_{\Bbb R}$, and put $g := \exp^{\Bbb C}
\mathcal{Y}\in T$, where the element $\exp (\mathcal{Y}/2)$ in $T$ 
is written as $\exp^{\Bbb C} \mathcal{Y}$ by abuse of terminology.
Recall that the $T$-action on $V_m$ 
is a right action, though it can be viewed also as a left action. 
Put $h_g := h\cdot g$ for simplicity.
Using the notation in
Definition 2.3, we write
$s_{k,i} = s_{l(k,i)}$,
$k =1,2,\dots,\nu_m$; 
$i = 1,2,\dots,n_k$. Then for a fixed $k$, 
 $\int_M | s_{k,i}|^{2}_{h^{m}_g}\,g^* \omega^n = |\chi_k
(\exp^{\Bbb C} \mathcal{Y})|^{-2}$ is
independent of the choice of $i$.
Put
$$ 
Z (q, \omega ; \mathcal{Y} )\;  := \;  \frac{n!}{m^n}\sum_{j=0}^{N_m}\; 
| s_j |^{\,2}_{h_g^m}\; =\; 
g^*\left \{ \frac{n!}{m^n}\sum_{k=1}^{\nu_m}\; |\chi_k (\exp^{\Bbb C} \mathcal{Y})|^{-2}\; 
\sum_{i=1}^{n_k} |s_{k,i}    |^2_{h^m}
\right \},
\quad \mathcal{Y}\in \frak t_{\Bbb R}.
$$

 For extremal K\"ahler manifolds, the following generalization of  \cite{M4} 
allows us to approximate arbitrarily 
some critical metrics relative to T:

\medskip
{\bf Theorem B.}\;  {\em  
  Let $\omega_0 = c_1(L; h_0 )$ be an extremal K\"ahler metric in the 
class $c_1(L)_{\Bbb R}$ with extremal K\"ahler vector field $\mathcal{V}$.  Then for some $T\in
\Delta$, 
there exist a sequence of vector fields 
$\mathcal{Y}_k\in \frak t_{\Bbb R}$, a formal power series $C_q$ in $q$ with real coefficients 
$(${\rm cf. Section} $6)$,
and smooth real-valued functions $\varphi_k$, $k = 1$,$2$,\dots,  on $M$ such that 
$$
Z  (q, \omega (\ell ); \mathcal{Y} (\ell )) 
\; =\; C_q + 0(q^{\ell +2}), 
\leqno{(2.8)}
$$
where
$\mathcal{Y} (\ell ) := (\sqrt{-1}\,\mathcal{V}/2)\,q^2\, +\, \Sigma_{k=1}^{\ell}\, q^{k+2}
\mathcal{Y}_k$, 
$h(\ell ):= h_0 \exp (-\Sigma_{k=1}^{\ell} q^k\varphi_k )$, and $\omega (\ell ) :=
c_1(L; h(\ell ))$.}

\medskip
The  equality (2.8) above means that there exists a positive real constant $A_{\ell}$ 
independent of $q$ such that 
$\| Z  (q, \omega (\ell ); \mathcal{Y} (\ell )) - C_q\|_{C^0(M)} \leq A_{\ell} q^{\ell + 2}$
for all $q$ with $0 \leq q \leq 1$. By \cite{Z1}, for every nonnegative integer $j$, 
a choice of a larger constant $A = A_{j,\ell} >0$ keeps Theorem B still valid even if
the $C^0(M)$-norm is replaced by the $C^j(M)$-norm.

\section{A stability criterion}

In this section, some stability criterion will be given as a preliminary.
In a forthcoming paper \cite{M5}, we actually use a stronger version of Theorem 3.2
which guarantees the stability only by checking  the closedness of orbits through a point
for special one-parameter
subgroups ``perpendicular'' to the isotropy subgroup.
Now, for a
connected reductive algebraic group $G$,  defined over
$\Bbb C$, we consider a representation of $G$ on an $N$-dimensional complex vector
space $W$.
We fix a maximal compact subgroup $G_c$ of $G$.
Moreover, 
let $\Bbb C^*$ be a one-dimensional algebraic
torus  with the maximal compact subgroup $\, S^1$.

\medskip
{\em Definition $3.1$.}
(a) An algebraic group homomorphism 
$\lambda :  \Bbb C^* \to G$
is said to be a {\it special one-parameter subgroup\/} of $G$,
if  the image
$\lambda (S^1)$ is contained in
$G_c$.  

\smallskip\noindent
 (b) A point $w\neq 0$ in $W$ is said to be
{\it stable\/}, if
the orbit $G \cdot w$ is closed in $W$. 

\medskip

Later, we apply the following stability criterion  to the case where
$W = W^*_m$ and $G = G_m$.
 Let $w\neq 0$ be a point in $W$.

 \medskip
{\bf Theorem 3.2.} 
{\em  A point $w$ as above is stable if and only if there exists a point $w'$ in the orbit 
$G\cdot w$ of $w$ such that $\lambda (\Bbb C^* )\cdot w'$
is closed in
$W$ for every special one-parameter subgroup
$\lambda : \Bbb C^* \to G$ of $G$.}
 
\medskip
{\rm Proof}. \,  We prove this by induction on $\dim (G\cdot w) $. If $\dim (G \cdot w) =
0$,  the statement of the above theorem is obviously true.
Hence, fixing a positive integer $k$, assume that the statement is true for all 
 $0 \neq w \in W$ such that
$\dim (G \cdot w) < k$.
Now, let $0 \neq w \in W$ be such that $\dim (G \cdot w) = k$, and 
the proof is reduced to showing the statement for such a point $w$.
Let $\Sigma (G)$ be the set of all special one-parameter subgroups of $G$. Fix
a $G_c$-invariant Hermitian metric $\|\; \|$ on $W$. 
The proof is divided into three steps:

\smallskip\noindent
{\em Step \/$1$}:
First, we prove  ``only if'' part of Theorem 3.2.
Assume that $w$ is stable. Since  $G\cdot w$ is closed in $W$,
 the nonnegative function on $G\cdot w $ defined by
$$
G\cdot w\, \owns\,  g \cdot w\; \mapsto \; \| g \cdot w \| \,\in \,\Bbb R, \qquad g \in G,
\leqno{(3.3)}
$$
has a critical point at some point $w'$ in $G \cdot w$. Let $\lambda\in \Sigma (G)$,
and it suffices to show the closedness of $\lambda (\Bbb C^* ) \cdot w'$ in $W$.
We may assume that $\dim \lambda (\Bbb C^* )\cdot w' >0$. 
Then by using the coordinate system associated to an  orthonormal basis for $W$, 
we can write $w'$ as $(w'_0, \dots , w'_r, 0, \dots, 0)$ in such a way that
$w'_{\alpha} \neq 0$ for all $0\leq \alpha \leq r$ and that
$$
\lambda (e^{ t}) \cdot w' \; = \; (e^{ t\gamma_0 }w'_0, 
\dots, e^{ t \gamma_r }w'_r, 0,\dots,
0),
\qquad t \in \Bbb C,
$$
where $\gamma_{\alpha}$, $\alpha = 0,1,\dots,r$, are integers independent of the choice of 
$t$ in $\Bbb C$. Since the closed orbit $G\cdot w$ does not contain the origin of $W$, 
the inclusion $\lambda (\Bbb C^* )\cdot w' \subset G\cdot w$ shows that $r \geq 1$
and that the coincidence $\gamma_0 =\gamma_1 = \dots =\gamma_r$ cannot occur.
In particular, 
$$
f (t):= \log \| \lambda (e^{t}) \cdot w'  \|^2 = \log \left ( e^{ 2t\gamma_0 } |w'_0|^2 + 
e^{ 2t\gamma_1 } |w'_1|^2 + \dots + e^{ 2t\gamma_r } |w'_r|^2 \right ),
\quad t \in \Bbb R,
$$
satisfies $f''(t) >0$ for all $t$. Moreover, since the function in (3.3) has a
critical point at $w'$, we have  $f'(0) =
0$. It now follows that
$\lim_{t\to +\infty} f(t) = + \infty$ and $\lim_{t\to -\infty} f(t) = + \infty$.
Hence $\lambda (\Bbb C^* )\cdot w' $ is closed in $W$, as required.

\smallskip\noindent
{\em Step \/$2$}:
To prove  ``if'' part of Theorem 3.2, we may assume that  $w = w'$ without loss of
generality. Hence, suppose that $\lambda (\Bbb C^* )\cdot w$ is closed in $W$
for every $\lambda \in \Sigma (G)$.
It then suffices to show that $G\cdot w$ is closed in $W$.
For contradiction, assume that $G\cdot w$ is not closed in $W$.
Since the closure of  $G\cdot w$ in $W$ always contains a closed orbit $O_1$ in $W$,
by $\dim O_1 < \dim (G\cdot w) = k$, the induction hypothesis shows that there exists a
point
$\hat{w}\in O_1$ such that
$$
\text{$\lambda (\Bbb C^* )\cdot \hat{w}$ is closed in $W$ for 
every $\lambda\in \Sigma (G)$.}
\leqno{(3.4)}
$$
Moreover, there exist elements $g_i $, $i = 1,2$, \dots, in $G$ such that
 $g_i \cdot w$ converges to $\hat{w}$ in $W$. 
Then for each $i$, we can write 
$g_i = \kappa'_i \cdot \exp (2\pi A_i )\cdot \kappa_i$ 
for some $\kappa_i$, $\kappa'_i \in G_c$ and 
for some $A_i \in \frak a$, where $2\pi \sqrt{-1}\,\frak a$ is the Lie algebra of
some maximal compact torus  in $G_c$.
Let $2\pi \sqrt{-1}\,\frak a_{\Bbb Z}$ be the kernel of the exponential map of the Lie algebra 
$2\pi \sqrt{-1}\,\frak a$, and
put $\frak a_{\Bbb Q} := \frak a_{\Bbb Z} \otimes \Bbb Q$.
 Replacing $\{\kappa_i\}$ by its 
subsequence if necessary,
we may assume that 
$$
\kappa_i \to \kappa_{\infty}
\; \text{ and }\;
\{\exp (2\pi A_i )\cdot \kappa_i \} \cdot w \to w_{\infty},
\qquad
\text{as $i \to \infty$},
\leqno{(3.5)}
$$
for some $\kappa_{\infty} \in G_c$ and $w_{\infty}\in G_c\cdot \hat{w}$. 
Then by (3.4), the orbit $\lambda (\Bbb C^* )\cdot w_{\infty}$ is also closed in $W$ for 
every 
$\lambda \in \Sigma ( G)$.  
Let $\frak a_{\infty}$ denote the Lie
subalgebra of $\frak a$ consisting of all elements in $\frak a$ whose associated vector fields
on $W$
 vanish at
$\kappa_{\infty}\cdot w$. For a Euclidean metric on $\frak a$ induced from a suitable bilinear from
on $\frak a_{\Bbb Q}$ defined over $\Bbb Q$, we write 
$\frak a$ as a direct sum $\frak a_{\infty}^{\perp}\oplus  \frak a^{}_{\infty}$, 
where $\frak a_{\infty}^{\perp}$ is the
orthogonal complement  of $\frak a^{}_{\infty}$ in $\frak a$.
 Let $\bar{A}_i$ be the image of $A_i$ under the orthogonal projection
$$
\operatorname{pr}_1 :\;\frak a\, (=   \frak a_{\infty}^{\perp}\oplus  \frak a^{}_{\infty})
\to\, \frak a^{\perp}_{\infty}, 
\qquad A \mapsto \bar{A}:= \operatorname{pr}_1(A).
$$
Note that $\{ \exp (2\pi A_i )\cdot \kappa_{\infty} \} \cdot w = 
\{ \exp (2\pi \bar{A}_i )\cdot \kappa_{\infty} \} \cdot w$.
Hence, 
\begin{align*}
 &\limsup_{i\to \infty} \| \, \exp\left\{ 2\pi 
{\operatorname{Ad}}(\kappa_{\infty}^{-1})\bar{A}_i \right\}\cdot
w\| \; =\; \limsup_{i\to \infty} \| \left \{ \exp (2\pi 
A_i )\cdot  \kappa_{\infty} \right\}\cdot
w\|     \tag{3.6}    \\
&\; \leq \; \lim_{i \to \infty}\| \left \{ \exp (2\pi A_i )\cdot \kappa_i  \right \} \cdot w \|\;
=\; \| w_{\infty} \| \;  < \; +\infty.
\end{align*}

\smallskip\noindent
{\em Step $3$}:
Since $\lambda (\Bbb C^* )\cdot w$ is closed in $W$ for 
every 
$\lambda\in \Sigma (G)$,  by the boundedness in (3.6),
$\{ \bar{A}_i\}$ is  a bounded sequence in 
$\frak a_{\infty}^{\perp}$ (see Remark 3.7 below).
Hence, for some element $A_{\infty}$ in $\frak a_{\infty}^{\perp}$,
replacing $\{\bar{A}_i\}$ by its subsequence if necessary, we may
assume that $\bar{A}_i \to A_{\infty}$ as $i\to\infty$.
Then by (3.5), 
$$
w_{\infty}\; =\; \lim_{i\to \infty}\; \{\exp (2\pi \bar{A}_i )\cdot \kappa_i \} \cdot w \; =\;
\{\exp (2\pi \bar{A}_{\infty} )\cdot \kappa_{\infty} \} \cdot w.
$$
Since  we have
$\exp (2\pi \bar{A}_{\infty} ) \in G$, the point
$w_{\infty}$ in $O_1$
belongs to the orbit $G\cdot w$. This contradicts 
 $O_1\cap (G\cdot w) = \emptyset$, as required.
The proof of Lemma 3.2 is now complete.

\medskip
{\em Remark $3.7$}. \;The
boundedness of the sequence $\{\bar{A}_i\}$ in $\frak a_{\infty}^{\perp}$ 
in Step 3 above can be seen as follows:
For contradiction, we assume that the sequence $\{\bar{A}_i\}$ is unbounded.
Put $v := \kappa_{\infty}\cdot w$ for simplicity.
Then by (3.6), we first observe that
$$
\limsup_{i\to\infty} \| \exp (2\pi \bar{A}_i )\cdot v \| < + \infty.
\leqno{(3.8)}
$$
Since $2\pi \sqrt{-1}\, \frak a_{\infty}$ is the Lie algebra of the isotropy
subgroup  of 
the compact torus $\exp (2\pi \sqrt{-1}\, \frak a )$ at $v$, 
both $\frak a_{\infty}$ and $\frak
a_{\infty}^{\perp}$  are defined over $\Bbb Q$ in $\frak a$.
By choosing a complex coordinate system of $W$, we can write $v$ as
$(v_0, \dots , v_r, 0, \dots, 0)$ for some integer $r$ with $0\leq r \leq \dim W-1$
such that $v_{\alpha} \neq 0$ for all $0\leq \alpha \leq r$ and
that
$$
\exp\,(2\pi \bar{A})\cdot v \; = \; (e^{2\pi \chi_0 (\bar{A}) }v_0,  \dots , 
e^{2\pi \chi_r (\bar{A})} v_r,
0,
\dots, 0 ),
\qquad \bar{A} \in \frak a^{\perp}_{\infty},
\leqno{(3.9)}
$$
where $\chi_{\alpha} : \frak a_{\infty}^{\perp} \to \Bbb R$, $\alpha = 0,1,\dots, r$, are 
additive characters defined over $\Bbb Q$. Put $n:= \dim_{\Bbb R} \frak a^{\perp}_{\infty}$,
and let $(\frak a^{\perp}_{\infty})_{\Bbb Q}$ denote the set of all rational points in 
$\frak a^{\perp}_{\infty}$.
Let us now identify
$$
\frak a^{\perp}_{\infty} \;=\; \Bbb R^n \qquad \text{ and }\qquad 
(\frak a^{\perp}_{\infty})_{\Bbb Q}\; =\; \Bbb Q^n,
$$
as vector spaces.
Since the orbit $\Bbb \lambda (\Bbb C^*)\cdot w$ is closed in $W$
for all special one-parameter subgroups $\lambda :\Bbb  C^* \to G$ of $G$, 
the same thing is true also for $\Bbb \lambda (\Bbb C^*)\cdot v$. Hence,
$$
\Bbb Q^n \setminus \{0\}\; \subset\; \bigcup_{\alpha, \beta =0}^r \; U_{\alpha\beta},
\leqno{(3.10)}
$$
where
$U_{\alpha\beta} :=\{\,A\in\frak a\,;\, 
\chi_{\alpha} (A) > 0 > \chi_{\beta}(A)\,\}$. 
Note that the boundaries of the open sets $U_{\alpha\beta}$, $1\leq \alpha\leq r$, 
$1\leq \beta\leq r$, in $\Bbb R^n$ sit in the union of
$\Bbb Q$-hyperplanes 
$$
H_{\alpha} \; :=\; \{\, \chi_{\alpha} = 0\,\},
\qquad \alpha =0,1,\dots, r,
$$
in $\Bbb R^r$. Since  an intersection of any finite number of 
hyperplanes 
$H_{\alpha}$, $\alpha =0,1,\dots,
r$, has dense rational points,  (3.10) above easily implies
$$
\Bbb R^n \setminus \{0\}\; =\; \bigcup_{\alpha, \beta =0}^r \; U_{\alpha\beta}.
\leqno{(3.11)}
$$
Replacing $\{\bar{A}_i \}$ by its suitable subsequence if necessary, 
we may assume that there exists an element $A_{\infty}$ in $\frak a_{\infty}^{\perp}\,
 (=\Bbb R^n )$ with $\| A_{\infty} \|_{\frak a} = 1$ such that 
$$
\lim_{i \to \infty} \, \frac{\bar{A}_i}{\;\| \bar{A}_i\|_{\frak a}}\; =\; A_{\infty},
$$
where $\|\;\|_{\frak a}$ denotes the Euclidean norm for $\frak a$ as in Step 2
in the proof of Theorem 3.2.
By (3.11), there exist $\alpha, \, \beta\in \{0,1,\dots, r\}$ such that
$A_{\infty} \in U_{\alpha\beta}$,
and in particular $\chi_{\alpha}(A_{\infty}) > 0$.
On the other hand, $\limsup_{i\to\infty} \| \bar{A}_i\|_{\frak a} \, =\, +\infty$ by our
assumption. 
Thus, 
$$
\limsup_{i\to\infty}\chi_{\alpha} (\bar{A}_i) \; =\; 
\limsup_{i \to\infty} \,\{\,\| \bar{A}_i\|_{\frak a}\cdot \chi_{\alpha} (\bar{A}_i/\|
\bar{A}_i\|_{\frak a})\,\}\; =\;  
(\limsup_{\i\to\infty}\| \bar{A}_i\|_{\frak a} )\,\chi_{\alpha} (A_{\infty})\;=\;+\infty,
$$
in contradiction to  (3.8) and (3.9),
as required.

\section{The Chow norm}

Take an algebraic torus $T \subset \operatorname{Aut}^0(M)$, and 
let $\iota : \operatorname{SL}(V_m) \to \operatorname{PGL}(V_m)$
be the natural projection, where we regard $\operatorname{Aut}^0(M)$
as a subgroup of $\operatorname{PGL}(V_m)$ via the Kodaira embedding
$\Phi_m : M \hookrightarrow \Bbb P^* (V_m)$, $m \gg 1$.
In this section,  
we fix  a $\tilde{T}_c$-invariant Hermitian metric $\rho $ on $V_m$,
where $\tilde{T}_c$ is the maximal compact subgroup of $\tilde{T}:= \iota^{-1}(T)$.
Obviously, in terms of this metric, $V(\chi_k ) \perp
V(\chi_{\ell})$ if $k \neq \ell$.
Using Deligne's pairings (cf. \cite{Deligne}, 8.3), 
Zhang (\cite{Z2}, 1.5) defined a special type of norm  on $W^*_m$,
called the {\it Chow norm\/},  as
a nonnegative real-valued function
$$
W^*_m \owns w \;\longmapsto \;\|w\|^{}_{\operatorname{CH}(\rho )}\in \Bbb R^{}_{\geq 0},
\leqno{(4.1)}
$$
with very significant properties described below. 
First, this is a norm, so that it  has the only zero at
the origin satisfying the homogeneity condition 
$$
\|c\, w\|^{}_{\operatorname{CH}(\rho )} 
= |c| \cdot \|w\|^{}_{\operatorname{CH}(\rho )}
\qquad \text{for all $(c, w) \in \Bbb C \times W^*_m$.}
$$
For the group $\operatorname{SL}(V_m)$,
we consider the maximal compact subgroup
$\operatorname{SU}(V_m; \rho)$.
For a special one-parameter subgroup
$$
\lambda : \Bbb C^* \to \operatorname{SL}(V_m)  
$$
of $\operatorname{SL}(V_m)$,
there exist integers $\gamma_j$, $j = 0,1,\dots, N_m$, and an orthonormal basis $\{s_0,
s_1, \dots , s_{N_m}\}$ for $(V_m, \rho )$ such that, for all $j$, 
$$
\lambda_z \cdot s_j = e^{z\gamma_j} s_j, 
\qquad z \in \Bbb C,
\leqno{(4.2)}
$$
where $\lambda_z := \lambda (e^z )$.  Recall that  the subvariety 
 $M_m$  in $\Bbb P^*(V_m)$ is the image of
the Kodaira embedding $\Phi_m : M \hookrightarrow  \Bbb P^*(V_m)$
defined by
$$
\Phi_m (p)  \; =\; (s_0(p):s_1(p):\dots : s_{N_m}(p)),
\qquad p \in M,
\leqno{(4.3)}
$$
where $\Bbb P^*(V_m)$ is identified with $\Bbb P^{N_m}(\Bbb C ) = \{(z_0:z_1:
\dots : z_{N_m})\}$.  
Put $M_{m,t} := \lambda_t (M_m)$ for each $t\in \Bbb R$.
As in Section 2, $\hat{M}_{m,t}:=\lambda_t\cdot \hat{M}_m $  is the nonzero point of $W^*_m$
sitting over the Chow point of the irreducible reduced cycle $M_{m,t}$ 
on $\Bbb P^*(V_m)$.
Then (cf. \cite{Z2},  1.4, 3.4.1)
$$
\frac{d}{dt}\left (\log \| \hat{M}_{m,t} \|^{}_{\operatorname{CH}(\rho )}\right ) \; =\; 
(n+1) \int_{M}   
\frac{\Sigma_{j=0}^{N_m}\, \gamma_j \,| \lambda_t \cdot s_j|^2} 
{\Sigma_{j=0}^{N_m}\, |\lambda_t\cdot s_j|^2} \,
(\Phi_m^*\lambda_t^*\omega_{\operatorname{FS}})^{n},  \leqno{(4.4)} 
$$
where $\omega_{\operatorname{FS}}$ is the Fubini-Study form 
$(\sqrt{-1}/2\pi)\partial\bar{\partial}\log (\Sigma_{j=0}^{N_m} |z_j|^2)$ on $\Bbb
P^*(V_m)$, and we regard $\lambda_t$ as a linear transformation of
$\Bbb P^*(V_m)$ 
induced by (4.2).  Note that the term $\Phi_m^*\lambda_t^*\omega_{\operatorname{FS}}$ 
above is just
$(\sqrt{-1}/2\pi )\partial\bar{\partial}\log (\Sigma_{j=0}^{N_m}\, |\lambda_t\cdot s_j|^2)$.
Put $\Gamma :=2\pi \sqrt{-1}\,\Bbb Z$. By setting
$$
\Bbb C /\Gamma \;=\; \{\,   t+ \sqrt{-1}\, \theta  \,;\, t \in \Bbb R, \; \theta
\in \Bbb R /(2\pi \Bbb Z)\,\},
$$
we consider the complexified situation.
 Let $\eta : M \times \Bbb C/\Gamma \to \Bbb P^*(V_m)$ be the map sending
each $(p, t+ \sqrt{-1}\,\theta )$ in $M \times \Bbb C/\Gamma$ to 
$ \lambda_{t+ \sqrt{-1}\,\theta} \cdot \Phi_m (p)$ 
in $\Bbb P^*(V_m)$.  For simplicity, we put
$$
Q \; := \; \frac{\Sigma_{j=0}^{N_m}\, \gamma_j e_{}^{2t\gamma_j}| s_j|^2} 
{\Sigma_{j=0}^{N_m}\, e_{}^{2t\gamma_j}| s_j|^2}
\; \left ( =\; \frac{\Sigma_{j=0}^{N_m}\, \gamma_j \,| \lambda_t \cdot s_j|^2} 
{\Sigma_{j=0}^{N_m}\, |\lambda_t\cdot s_j|^2} \right ).
$$
We further put $ z := t+ \sqrt{-1}\,\theta$. 
For the time being, on the total complex manifold
$M \times \Bbb C/\Gamma$, the $\partial$-operator and the $\bar{\partial}$-operator
will be written simply as  $\partial$ and $\bar{\partial}$ respectively,
while on $M$, they will be denoted by $\partial_M$ and $\bar{\partial}_M$ respectively. Then
$$
\eta^*\omega_{\operatorname{FS}} \; 
=\;\Phi_m^* \lambda_t^*\omega_{\operatorname{FS}}\, 
+\, \frac{\sqrt{-1}}{2\pi}\,(\partial^{}_M Q \wedge d\bar{z}+dz\wedge \bar{\partial}^{}_MQ
)\, +\, \frac{\sqrt{-1}}{4 \pi}\, \frac{\partial Q}{\partial t}\,dz\wedge d\bar{z}.
$$
For $0\neq r\in \Bbb R$, we consider the $1$-chain $I_{r} := [0, r ]$,
where  $[0, r ]$ means the $1$-chain $\, -[r, 0]$ if $r <0$.
Let $\operatorname{pr}: \Bbb C/\Gamma\to \Bbb R$
be the mapping sending each $t + \sqrt{-1}\,\theta$ to $t$.
We now put $B_{r} : = \operatorname{pr}^*I_{r}$.
Then $\int_{M\times B_{r}} \eta^*\omega_{\operatorname{FS}}^{n+1}$ is nothing but
\begin{align*}
& (n+1)\int_{0}^{r} dt\int_M
 \left (\, \frac{\partial Q}{\partial t}\,\Phi_m^*
\lambda_t^*\omega^n_{\operatorname{FS}}\, +\,
\frac{\sqrt{-1}}{\pi}\,\bar{\partial}^{}_M Q \wedge
\partial^{}_M Q\wedge \,n\,\Phi_m^* \lambda_t^*\omega^{n-1}_{\operatorname{FS}}
\right )\\
&=\; \int_0^{r} \frac{d^2}{dt^2}\left (\log \| \hat{M}_{m,t}
\|^{}_{\operatorname{CH}(\rho )}\right ) dt\; =\; \frac{d}{dt}\left (\log \| \hat{M}_{m,t}
\|^{}_{\operatorname{CH}(\rho )}\right )\; \big |^{t=r}_{t=0},
\end{align*}
and by assuming $r \geq 0$, we obtain the following convexity formula:
$$
\frac{d}{dt}\left (\log \| \hat{M}_{m,t}
\|^{}_{\operatorname{CH}(\rho )}\right )\; \big |^{t=r}_{t=0} \; =\; \int_{M\times B_{r}}
\eta^*\omega_{\operatorname{FS}}^{n+1} \;\geq \; 0.
\leqno{\bf \quad Theorem \; 4.5.} 
$$

{\em Remark\, $4.6$}.  Besides special one-parameter subgroups of $\operatorname{SL}(V_m)$,
we also consider a little more general smooth path
$\lambda_t$, $t \in \Bbb R$, 
in  $\operatorname{GL}(V_m)$ written explicitly  by
$$
\lambda_t \cdot s_j = e^{t\gamma_j + \delta_j} s_j, 
\qquad j = 0,1, \dots, N_m,
$$
where $\gamma_j$, $\delta_j \in \Bbb R$ are not necessarily rational.
In this case also, we easily see that the formula (4.4) and Theorem 4.5 are still valid.
 
\section{Proof of Theorem A}

The statement of Theorem A is divided into  ``if'' part, ``only if'' part, and the uniqueness
part. We shall prove these three parts separately.

\medskip
Proof of ``{\it if}'' part. \; Let $\omega\in \mathcal{S}$ be a critical metric relative to $T$.
Then by Definition 2.6, in terms of the Hermitian metric defined in (2.4), 
there exists an admissible normal basis $\{s_0, s_1, \dots, s_{N_m}\}$ for $V_m$ 
of index $b$ such that 
the associated function $E_{\omega, b}$ has a constant value $C$ on $M$.
By operating
$(\sqrt{-1}/2\pi) \partial \bar{\partial}\log$ on the identity $E_{\omega, b} = C$, we have
$$
\Phi_m^* \omega_{\operatorname{FS}} \; =\; m\, \omega.
\leqno{(5.1)}
$$
Besides the Hermitian metric defined in (2.4), we shall now define another Hermitian metric on
$V_m$. By the identification  $V_m\cong \Bbb C^{N_m}$ via the basis $\{s_0, s_1, \dots, s_{N_m}\}$,
the standard Hermitian metric on $\Bbb C^{N_m}$ induces a Hermitian metric $\rho$ 
on $V_m$.
As a maximal compact subgroup of $G_m$, we choose $(G_m)_c$ as in Section 2
by using the metric defined in (2.4). Then the Hermitian metric $\rho$ is also 
preserved by the $(G_m)_c$-action on $V_m$.   Let 
$$
\lambda : \Bbb C^* \to G_m
$$ 
be a special one-parameter subgroup of $G_m$.
By the  notation $l (k, i)$ as in Definition 2.3,
we put $s_{k,i} := s_{l(k,i)}$.
If necessary, replacing $\{s_0, s_1, \dots, s_{N_m}\}$ by another
admissible normal basis for $V_m$ of the same index $b$,
 we may assume without loss of generality  that
there exist integers $\gamma_{k,i}$, $i = 1,2,\dots, n_k$, satisfying
$$
\lambda_t\cdot s_{k,i} \; =\; e_{}^{t\gamma_{k,i}} s_{k,i},
\qquad t \in \Bbb C,
\leqno{(5.2)}
$$
where $\lambda_t := \lambda (e^t)$ is as in (4.2), and 
the equality $\Sigma_{i=1}^{n_k} \gamma_{k,i} = 0$ is required to hold for every $k$.
Put $\gamma_{k, i} = \gamma_{ l (k, i)}$ for simplicity.
Then by (4.4) and (5.1),
\begin{align*}
&\frac{d}{dt}\left ( \log \| \hat{M}_{m,t}\|_{\operatorname{CH}(\rho )}\right )_{|t = 0}
\; =\; (n+1) \int_M \frac{\Sigma_{j=0}^{N_m} \gamma_j |s_j |^2}{\Sigma_{j=0}^{N_m} 
|s_j |^2}(\Phi_m^*  \omega_{\operatorname{FS}})^n\\
&= \; (n+1)\,m^n \int_M \frac{\Sigma_{j=0}^{N_m} \gamma_j  |s_j
|_{h^m}^2}{\Sigma_{j=0}^{N_m} 
|s_j |_{h^m}^2}\,\omega^n
\; =\;
(n+1)\,m^n\int_M\frac{\Sigma_{k=1}^{\nu_m}(\Sigma_{i=1}^{n_k}\gamma_{k,i} |s_i |^2_{h^m}) }
{ E_{\omega, b}} 
\omega^n\\
&= \; \frac{(n+1)\, m^n}{C}\int_M 
\Sigma_{k=1}^{\nu_m}(\Sigma_{i=1}^{n_k}\gamma_{k,i} |s_i |^2_{h^m})
\,\omega^n \; =\; \frac{(n+1)\, m^n}{C}\,
\Sigma_{k=1}^{\nu_m}\, b_k(\Sigma_{i=1}^{n_k}\gamma_{k,i}) \; =\; 0.
\end{align*}
Note also that, by Theorem 4.5, we have  $c := (d^2/dt^2) ( \log \|
\hat{M}_{m,t}\|_{\operatorname{CH}(\rho )} )_{|t =0}\, \geq \, 0$.

\smallskip\noindent
{\em Case $1$}:\; 
If $c$ is positive, then 
$\lim_{t\to -\infty} \| \hat{M}_{m,t}\|_{\operatorname{CH}(\rho )} = +\infty
= \lim_{t\to +\infty} \| \hat{M}_{m,t}\|_{\operatorname{CH}(\rho )}$,
and in particular $\lambda (\Bbb C^*)\cdot \hat{M}_m$ is closed.

\smallskip\noindent
{\em Case $2$}: \;
If $c$ is zero, then by applying
Theorem 4.5 infinitesimally, we see that 
$\lambda (\Bbb C^*)$ preserves the subvariety $M_m$ in $\Bbb
P^* (V_m)$, and moreover by 
$(d/dt ) ( \log \| \hat{M}_{m,t}\|_{\operatorname{CH}(\rho )})_{|t = 0} = 0$, 
the isotropy representation of $\lambda (\Bbb C^*)$ on the complex line $\Bbb C \hat{M}_m$ 
is trivial. Hence, $\lambda (\Bbb C^*)\cdot \hat{M}_m$ is a single point, 
and in particular closed.

\smallskip
Thus, these two cases together with
Theorem 3.2 show that the subvariety $M_m$ of $\Bbb P^* (V_m)$ is stable relative to $T$,
as required.

\medskip
{\em Remark\, $5.3$}. About the one-parameter subgroup 
$\{\lambda_t\,;\, t \in \Bbb R\}$ of $G_m$, we
consider a more general situation that
$\gamma_{k,i}$ in (5.2) are just real numbers
which are not necessarily rational.
The above computation together with Remark 4.6 shows that, even in this case,
$(d/dt)_{t=0} ( \log \| \hat{M}_{m,t}\|_{\operatorname{CH}(\rho )})$ vanishes.

\medskip
Proof of ``{\it only if}'' part. \;  Assume that the subvariety $M_m$ in $\Bbb P^* (V_m )$ is
stable relative to $T$. Take a  Hermitian metric $\rho$ for $V_m$ such that $V(\chi_k )\perp
V(\chi_{\ell})$ for  $k \neq \ell$. For this $\rho$, we consider the associated Chow norm.
Since the
orbit $G_m\cdot\hat{M}_m$ is closed in
$W_m$, the Chow norm restricted to this orbit attains an abosolute minimum.
 Hence, for some $g_0\in G_m$,
$$
0\; \neq\;  \| g_0\cdot \hat{M}_m\|^{}_{\operatorname{CH}(\rho )} \; \leq \;
\| g\cdot \hat{M}_m\|_{\operatorname{CH}(\rho )},
\qquad \text{ for all $g \in G_m$.}
$$
By choosing an admissible normal basis 
$\{ s_0, s_1, \dots, s_{N_m}\}$ for $(V_m; \rho )$ of index $(1,1,\dots,1)$,
we identify $V_m$ with  $\Bbb C^{N_m} =\{\,(z_0, z_1, \dots, z_{N_m})\,\}$.
Then $\operatorname{SL}(V_m)$ is identified with
$\operatorname{SL}(N_m +1;\Bbb C )$.
Let $\frak g_m$ be 
the Lie subalgebra of $\frak{sl}(N_m +1; \Bbb C)$
associated to the Lie subgroup 
$G_m$ of
$\operatorname{SL}(N_m +1;\Bbb C )$. We can now write $\,g_0 \, =\, \kappa'\cdot \exp
\{\operatorname{Ad}(\kappa ) D \} $ for some
$\kappa$, $\kappa' \in G_{m,c}$ and a  real diagonal matrix $D$
in $\frak g_m$.  By $\| \exp \{\operatorname{Ad}(\kappa ) D\} \cdot
\hat{M}_m\|_{\operatorname{CH(\rho )}}
 = \| g_0\cdot \hat{M}_m\|^{}_{\operatorname{CH}(\rho )}$, we have
$$
\;
\| \exp \{\operatorname{Ad}(\kappa ) D\}\cdot
\hat{M}_m\|_{\operatorname{CH(\rho )}}
\, \leq \, \| \exp \{ t \operatorname{Ad}(\kappa )A\}
\cdot \exp \{\operatorname{Ad}(\kappa )
D\}\cdot
\hat{M}_m\|_{\operatorname{CH(\rho )}},
\;\; t \in \Bbb R,
\leqno{(5.4)}
$$
for every real diagonal matrix $A$ in $\frak g_m$. 
For  $j =0,1,\dots, N_m$, we write the $j$-th diagonal element
  of $A$ and $D$ above as  $a_j$ and $d_j$,
respectively. Put $c_j := \exp\, d_j$ and $s'_j := \kappa^{-1}\cdot s_j$. Then 
$\{ s'_0, s'_1, \dots, s'_{N_m}\}$ is again an admissible normal basis for $(V_m,\rho )$ 
of index $(1,1,\dots,1)$.
By the notation in Definition 2.3,  we rewrite $s'_j$, $a_j$, $c_j$, $z_j$  as 
$s'_{k,i}$, $a_{k,i}$, $c_{k,i}$, $z_{k,i}$  by
$$
s'_{k,i} := s'_{l (k,i)},
\quad a_{k,i} := a_{l (k,i)},
\quad c_{k,i} := c_{l (k,i)},
\quad z_{k,i} := z_{l (k,i)},
$$
where $k = 1,2,\dots, \nu_m$ and $i = 1,2,\dots, n_k$.
By (5.4), the derivative at $t = 0$ of the right-hand side of (5.4) vanishes.
Hence by (4.4) together with Remark 4.6, 
fixing an arbitrary real diagonal matrix $A$ in $\frak g_m$, we have
$$
\int_M \frac{\Sigma_{k=1}^{\nu_m} \Sigma_{i=1}^{n_k} a_{k,i} c^{\,2}_{k,i}|s'_{k,i}|^2}
{\Sigma_{k=1}^{\nu_m} \Sigma_{i=1}^{n_k} c^{\,2}_{k,i}|s'_{k,i}|^2}\;
\Phi_m^*(\Theta^n )\; =\; 0
\leqno{(5.5)}
$$
where we set $\Theta :=  (\sqrt{-1}/2\pi )\partial\bar{\partial}\log 
(\Sigma_{k=1}^{\nu_m} \Sigma_{i=1}^{n_k} c^{\,2}_{k,i}|z_{k,i} |^2)$.
Let $k_0\in \{1,2,\dots, \nu_m\}$
and let $i_1$, $i_2 \in\{1,2,\dots, n_k\}$ with $i_1 \neq i_2$. Using
Kronecker's delta, we specify the real diagonal matrix  $A$  by
setting
$$
a_{k, i } \; =\; \delta_{k k_0}(\delta_{ii_1} - \delta_{i i_2}),\qquad 
k = 1,2,\dots, \nu_m; \; i =1,2,\dots, n_k.
$$
Apply (5.5) to this $A$, and let $(i_1, i_2)$ run through the set of all pairs of two distinct
elements in $\{1,2,\dots, n_k\}$. Then there exists a positive constant $b_k >0$ independent
of the choice of $i$ in $\{1,2,\dots, n_k\}$ such that
$$
\frac{N_m + 1}{m^n c_1(L)^n [M]}\int_M \frac{ c^{\,2}_{k,i}|s'_{k,i}|^2}
{\Sigma_{k=1}^{\nu_m} \Sigma_{i=1}^{n_k} c^{\,2}_{k,i}|s'_{k,i}|^2}\;
\Phi_m^*(\Theta^n )\; =\; b_k,
\qquad k = 1,2,\dots, \nu_m.
\leqno{(5.6)}
$$
The following identity (5.7) allows us to define (cf. \cite{Z2})
 a Hermitian metric $h_{\operatorname{FS}}$ on $L^m$ by
$$
|s |^2_{h_{\operatorname{FS}} } :=\;  \frac{(N_m + 1)}{c_1(L)^n[M]}\,
\frac{\Sigma_{k=1}^{\nu_m} \Sigma_{i=1}^{n_k}
\, |\, (s, s'_{k,i})_{\rho}|^2\,|s'_{k,i} |^2 }
{\Sigma_{k=1}^{\nu_m} \Sigma_{i=1}^{n_k} c^{\,2}_{k,i}|s'_{k,i}|^2},
\qquad s \in V_m.
\leqno{(5.7)}
$$
Then for this Hermitian metric,
it is easily seen that
$$
\Sigma_{j=0}^{N_m} | c_j s'_j |^2_{h_{\operatorname{FS}} }
= \Sigma_{k=1}^{\nu_m} \Sigma_{i=1}^{n_k} | c_{k,i} s'_{k,i} |^2_{h_{\operatorname{FS}} }
= (N_m +1)/c_1(L)^n[M].
\leqno{(5.8)}
$$
By operating $ (\sqrt{-1}/2\pi )\partial\bar{\partial}\log$ on  
both sides of (5.8), we obtain
$\Phi_m^* \Theta \; =\; c_1 (L^m; h_{\operatorname{FS}})$.
We now set $h := (h_{\operatorname{FS}})^{1/m}$ and $\omega := c_1 (L; h)$. Then
$$
\omega = (1/m)\, \Phi_m^*\Theta.
$$
Put $s''_{k,i} := c_{k,i}s'_{k,i}$, and  as in Definition 2.3, 
we write $s''_{k,i}$ as $s''_{l (k,i)}$.  Then by (5.8), we have the equality
$\Sigma_{j=0}^{N_m} | s_j'' |_{h^m}^2 = (N_m +1)/c_1(L)^n[M]$.  Moreover, in terms of the
Hermitian metric defined in (2.4), the equality (5.6) is interpreted as 
$$
\| s''_{k,i}\|^{\,2}_{L^2} = b_k, \qquad k = 1,2,\dots, \nu_m; \; i = 1,2,\dots, n_k,
$$
while by this together with  (5.8) above, we obtain
$\Sigma_{k=1}^{\nu_m}\, n_k b_k = N_m +1$, as required.

\medskip\noindent
Proof of {\it uniqueness}. \;  Let $\omega = c_1 (L; h)$ and $\omega' = c_1 (L; h')$ 
be critical metrics relative to $T$, and let 
$\{ s_j\,;\, j = 0,1,\dots, N_m\,\}$
and $\{ s'_j\,;\, j =0,1,\dots, N_m\,\}$ be respectively the 
associated admissible normal bases for $V_m$ of index $b$.
We use the notation in Definition 2.3. Then
$$
E_{\omega, b} := \Sigma_{k=1}^{\nu_m}\Sigma_{i=1}^{n_k}\, |s_{k,i} |^2_{h^m}
\quad\text{and}\quad
E_{\omega', b} := \Sigma_{k=1}^{\nu_m}\Sigma_{i=1}^{n_k}\, |s'_{k,i}|^2_{{h'}^m}
$$
take the same constant value $C := (N_m +1)/c_1(L)^n[M]$ on
$M$.  Note here that, by operating $(\sqrt{-1}/2\pi ) \partial \bar{\partial}\log$ 
on both of these identities,
we obtain
$$
m\omega = (\sqrt{-1}/2\pi )\partial\bar{\partial}\log
(\Sigma_{k=1}^{\nu_m}\Sigma_{i=1}^{n_k}\, |s_{k,i} |^2)
\;\; \text{and} \;\;
m\omega' = (\sqrt{-1}/2\pi )\partial\bar{\partial}\log
(\Sigma_{k=1}^{\nu_m}\Sigma_{i=1}^{n_k}\, |s'_{k,i} |^2).
$$
If necessary, we replace  each $s_{k, i}$
by $\zeta_k  s_{k,i}$ for a suitable complex number $\zeta_k$,
independent of  $i$, of absolute value 1.
Then for each $k = 1,2,\dots, \nu_m$, we may assume that there exist
a matrix $g^{(k)} = (g^{(k)}_{i\,\hat{i}}) \in \operatorname{GL}(n_k;\Bbb C)$ 
satisfying
$$
s'_{k,\hat{i}} \; =\;   \sum_{i=1}^{n_k} \, s_{k,i}\,g^{(k)}_{i\,\hat{i}},
$$
where $i$ and $\hat{i}$ always run through the integers in $\{1,2,\dots,n_k\}$.
Then the matrix $g^{(k)}$ above is written as $\kappa^{(k)} \cdot (\exp A^{(k)})\cdot 
(\kappa'{}^{(k)})_{}^{-1}$ 
for some real 
diagonal matrix $A^{(k)}$ and  
$$
\kappa^{(k)} = (\kappa^{(k)}_{i\,\hat{i}})\quad 
\text{and}\quad \kappa'{}^{(k)} = (\kappa'{}^{(k)}_{i\,\hat{i}})
$$
in $\operatorname{SU}(n_k )$. Let $a^{(k)}_i$ be the $i$-th diagonal element of $A^{(k)}$.
For each $\hat{i}$, we put $\tilde{s}_{k,\hat{i}} := \Sigma_{i=1}^{n_k} s_{k,i}\,\kappa_{i\,\hat{i}}^{(k)}$
and 
$\tilde{s}'_{k,\hat{i}} := \Sigma_{i=1}^{n_k} s'_{k,i}\,\kappa'_{i\,\hat{i}}{}^{(k)}$.
If necessary, we replace the bases $\{s_{k,1}, s_{k,2}, \dots, s_{k,n_k}\}$ and $\{s'_{k,1}, s'_{k,2}, \dots, s'_{k,n_k}\}$  
for $V (\chi_k )$   by the bases
$\{\tilde{s}_{k,1}, \tilde{s}_{k,2}, \dots, \tilde{s}_{k,n_k}\}$ and
 $\{\tilde{s}'_{k,1}, \tilde{s}'_{k,2}, \dots, \tilde{s}'_{k,n_k}\}$, respectively.
Then we may assume, from the beginning, that
$$
s'_{k,i} \; =\;  \{\exp\, a_i^{(k)}\}\, s_{k,i},
\qquad i =1,2,\dots, n_k.
$$
We now set $\tau_{k,i} :=s_{k,i}/\sqrt{b_k}$, 
and the Hermitian metric for $V_m$ defined in (2.4) will be denoted by $\rho$.
Then $\{\tau_{k, i}\,;\,k =1,2,\dots, \nu_m, i = 1.2.\dots, n_k\}$
is an admissible normal basis
of index $(1,1,\dots, 1)$   for $(V_m, \rho )$. 
Let $\{\lambda_t \, ;\, t\in \Bbb C\}$ be the smooth one-parameter family of elements
in $\operatorname{GL}(V_m)$ defined by
$$
\lambda_t \cdot \tau_{k,i} \; =\;  \{\exp (t\, a_i^{(k)})\}\sqrt{b_k}\,\tau_{k,i},
\qquad k = 1,2,\dots, \nu_m; \; i = 1,2,\dots, n_k.
$$
Put $\hat{M}_{m,t} := \lambda_t \cdot \hat{M}_m$, $0\leq t \leq 1$. Then by Remark 4.6
applied to the formula (4.4),
the derivative $\frak d (t) := (d/dt)( \log \| \hat{M}_{m,t}\|_{\operatorname{CH}(\rho )} )
/(n+1)$
at 
$t\in [0,1]$
is expressible as
$$
 \int_M \frac{\Sigma_{k=1}^{\nu_m}\Sigma_{i=1}^{n_k} a_i^{(k)} | \lambda_t 
\cdot \tau_{k,i} |^2}
{\Sigma_{k=1}^{\nu_m}\Sigma_{i=1}^{n_k}  | \lambda_t \cdot \tau_{k,i} |^2}
\left \{ (\sqrt{-1}/2\pi ) \partial\bar{\partial}
\log (\Sigma_{k=1}^{\nu_m}\Sigma_{i=1}^{n_k} |\lambda_t \cdot \tau_{k,i} |^2)\right \}^n
$$
Hence at $t = 0$, we see that 
$$
\frak d (0) = \int_M \Sigma_{k=1}^{\nu_m}\Sigma_{i=1}^{n_k} \{a_i^{(k)} |
s_{k,i} |_{h^m}^2/C\}(m\omega )^n = (m^n/C)  \, \Sigma_{k=1}^{\nu_m}
\{b_k\Sigma_{i=1}^{n_k} a_i^{(k)}\},
$$
while at $t = 1$ also, we obtain
$$
\frak d (1) =\int_M \Sigma_{k=1}^{\nu_m}\Sigma_{i=1}^{n_k} \{a_i^{(k)} |
s'_{k,i} |_{h'{}^m}^2/C\}(m\omega' )^n = (m^n/C)  \, \Sigma_{k=1}^{\nu_m}
\{b_k\Sigma_{i=1}^{n_k} a_i^{(k)}\}.
$$
Thus, $\frak d (0)$ coincides with $\frak d (1)$, while by Remark 4.6, we see from Theorem 4.5 
that
$(d^2/dt^2) \{\log \| \hat{M}_{m,t} \|^{}_{\operatorname{CH}(\rho )}\} \geq 0$ on $[0,1]$.
Hence, for all $t \in [0,1]$,
$$
\frac{d^2}{dt^2} \{\log \| \hat{M}_{m,t} \|^{}_{\operatorname{CH}(\rho )}\} = 0,
\qquad \text{ on $M$.}
$$
By Remark 4.6, the formula in  Theorem 4.5 shows that
$\lambda_t$, $t \in [0,1]$, belong to $H$ up to a positive scalar multiple.
Since $\lambda_1$ commutes with $T$, 
the uniqueness follows, as required.

\section{Proof of Theorem B}

Throughout this section, we assume that the first Chern class $c_1(L)_{\Bbb R}$ admits
an extremal K\"ahler metric $\omega_0 = c_1(L; h_0)$. Then  by a theorem of Calabi
\cite{Ca}, the identity component $K$ of the group of isometries of $(M, \omega_0)$ is a 
maximal compact connected subgroup of $H$, and we obtain $\omega_0 \in \mathcal{S}_K$
by the notation in the introduction.

\medskip
{\em Definition \/$6.1.$} \; 
For a $K$-invariant K\"ahler metric $\omega \in \mathcal{S}_K$
 on $M$ in the class $c_1(L)_{\Bbb R}$, we choose a Hermitian metric $h$ on $L$
such that
$\omega\, =\,c_1(L;h)$.
Then the power series  in $q$ given by the  right-hand side of (2.8)
will be denoted by $\Psi (\omega, q)$. Given $\omega$ and $q$, 
the power series  $\Psi (\omega, q)$ is independent of the choice of $h$.

\medskip
Let $\mathcal{D}_0$ be the Lichn\'erowicz operator as
defined in
\cite{Ca}, (2.1), for the  extremal K\"ahler manifold $(M, \omega_0 )$. 
Then by $\mathcal{V}\in \frak k$, the operator
$\mathcal{D}_0$  preserves the space $\mathcal{F}$ of all
real-valued smooth
$K$-invariant functions $\varphi$  such that $\int_M \varphi \omega_0^{\,n} = 0$.
Hence, we regard $\mathcal{D}_0$
just as an operator $\mathcal{D}_0 :
\mathcal{F} \to \mathcal{F}$, and the kernel in 
$\mathcal{F}$ of this restricted
operator will be denoted simply by $\operatorname{Ker} \mathcal{D}_0$.  
Then $\operatorname{Ker}
\mathcal{D}_0$  is a subspace of $\mathcal{K}_{\omega_0}$,
and we have an isomorphism
$$
e_{0} : \operatorname{Ker} \mathcal{D}_0 \; \cong \; \frak z,
\qquad \varphi \leftrightarrow 
e_{0} (\varphi ) := \operatorname{grad}_{\omega_0}^{\Bbb C}\varphi.
\leqno{(6.2)}
$$
 By the inner product $(\;,\;)_{\omega_0}$ defined in the introduction,
we write $\mathcal{F}$ as an orthogonal  direct sum $\operatorname{Ker} \mathcal{D}_{0} \oplus
\operatorname{Ker} \mathcal{D}_{0}^{\perp}$. 
We then consider the
orthogonal projection
$$
P\;  : \; \mathcal{F} \; (=\operatorname{Ker} \mathcal{D}_0 \oplus
\operatorname{Ker} \mathcal{D}_0^{\perp})
\;\to\;
\operatorname{Ker} \mathcal{D}_0.
$$
Now, starting from $\omega (0):= \omega_0$, we inductively define  a Hermitian
metric  $h(k)$,
a K\"ahler metric 
$\omega (k):= c_1(L; h(k))\in \mathcal{S}_K$, 
and  a vector field $\mathcal{Y}(k)\in
\sqrt{-1}\,\frak z$, $k =1,2,\dots$, 
 by 
\begin{equation}
\begin{cases}
&h(k)\;  :=\;  h(k-1) \exp (- q^k\varphi_k ),   \\
&\omega (k) \; =\; \omega (k-1)\, +\, (\sqrt{-1}/2\pi )
\,q^k\,\partial\bar{\partial}\varphi_k , 
 \\
&\mathcal{Y}(k) \; =\; \mathcal{Y}(k-1) \,+\,  \sqrt{-1}\; q^{k+2} e_{0} (\zeta_k ),
\end{cases} \tag{6.3}
\end{equation}
for appropriate  $\varphi_k \in \operatorname{Ker} \mathcal{D}_0^{\perp}$
and $\zeta_k \in\operatorname{Ker} \mathcal{D}_0$,
where   $\omega (k)$ and $\mathcal{Y}(k)$ are required to satisfy 
the condition (2.8) with $\ell$ replaced by $k$. We now set 
$g(k) := \exp^{\Bbb C} \mathcal{Y}(k)$.
Then
\begin{align*}
&\{ h(k)\cdot g(k)\}^{-m}h(k)^m\,\{ Z(q, \omega (k); \mathcal{Y}(k)) - C_{q}\}\\
&=\;\frac{n!}{m^n}\{ \Sigma_{j=0}^{N_m} |s_j |_{h(k)^m}\}
-C_{q}\{g(k)\cdot h(k)^{-m}\}h(k)^m\\
&= \; \Psi (\omega (k), q )\, -\,  C_{q}\, h(k)^m\,\{ (\exp^{\Bbb C} \mathcal{Y}(k)) \cdot h(k)^{-m}\},\\
& =\; \Psi (\omega (k), q )\, -\,  C_{q}\,\left \{ 1 +\, h(k)\, (\mathcal{Y}(k)/q) \cdot h(k)^{-1}\, +\,
 R(\mathcal{Y}(k); h(k))\right \},
\end{align*}
where $C_{q} = 1 + \Sigma_{k=0}^{\infty}\, \alpha_k q^{k+1}$ is a power series in $q$ with 
real coefficients $\alpha_k$ specified later,
and the last term $R(\mathcal{Y}(k); h(k)) := h(k)^m \Sigma_{j=2}^{\infty}\, 
\{ \mathcal{Y}(k)^j/j! \}\cdot h(k)^{-m}$ will be taken care of as a higher order term in $q$.
Consider the truncated term $C_{q,\ell} = 1 +
\Sigma_{k=0}^{\ell}\,
 \alpha_k q^{k+1}$.
Put
$$
\Xi (\omega (k), \mathcal{Y}(k), C_{q,k}) : = 
\Psi (\omega (k), q )\, -\,  C_{q, k}\,\left \{ \,1 -\,  (\mathcal{Y}(k)/q) \cdot \log h(k)\, +\,
 R(\mathcal{Y}(k); h(k))\,\right \}
$$ 
for each $k$. Then, in terms of $\omega (k)$, $\mathcal{Y}(k)$ and 
$C_{q,k}$, the condition (2.8) with $\ell$ replaced by $k$ is just
the equivalence 
$$
\Xi (\omega (k), \mathcal{Y}(k), C_{q,k}) \,\equiv \, 0, 
\qquad \text{modulo $q^{k+2}$}.
\leqno{(6.4)}
$$

\smallskip
We shall now  define $\omega (k)$, $\mathcal{Y}(k)$ and 
$C_{q,k}$  inductively in such a way that the condition (6.4) is satisfied.
If $k =0$, then we set $\omega (0) = \omega_0$, 
$\mathcal{Y}(0) = \sqrt{-1}\,q^2\mathcal{V}/2$
and $C_{q,0} =1+ \alpha_0 q$, where we put
$\alpha_0 := \{2c_1(L)^n[M]\}^{-1}\{\int_M \sigma_{\omega} \omega^n + 2 \pi F(\mathcal{V})\}$
for 
$\omega\in \mathcal{S}_K$. This $\alpha_0$
is obviously independent of the choice of $\omega$ in $\mathcal{S}_K$.
Then, modulo $q^2$, 
\begin{align*}
&\Psi (\omega (k), q )\, -\,  C_{q,0}\,\left \{ 1 -\,  (\mathcal{Y}(0)/q) \cdot \log h(0)\, +\,
 R(\mathcal{Y}(0); h(0))\right \}\\
& \equiv \; \left (1\, +\, \frac{\sigma_{\omega_0}}{2}\,q\right )
\,-\, (1 + \alpha_0 q) \left \{ 1 \,-\,  \,q\,h^{-1}_0
\sqrt{-1}\, (\mathcal{V}/2)\cdot h_0  \right \}\\
&\equiv \;
 \left (1\, +\, \frac{\sigma_{\omega_0}}{2}\,q\right )
\,-\, (1 + \alpha_0 q) \left \{ 1 + 
\left (\frac{\sigma_{\omega_0}}{2} -\alpha_0 \right )q \right \}
\;\equiv \; 0,
\end{align*}
and we see that  (6.4) is true for 
$k = 0$.  Here, the equality $h^{-1}_0
\sqrt{-1}\, (\mathcal{V}/2)\cdot h_0
= \alpha_0 - (\sigma_{\omega_0}/2)$ follows from a routine computation
(see for instance \cite{M1}).

\medskip
Hence, let $\ell \geq 1$ and assume (6.4)  for $ k
=\ell -1$. It then suffices to find $\varphi_\ell$, $\zeta_\ell$ and $\alpha_{\ell}$ 
satisfying  (6.4) for
$k = \ell$. 
 Put $\mathcal{Y}_{\ell} := \sqrt{-1}\,
e_{0}(\zeta_\ell )$. For each $(\varphi_\ell,  \zeta_\ell, \alpha_\ell)
\in \operatorname{Ker} \mathcal{D}_0^{\perp}\times \operatorname{Ker} \mathcal{D}_0
\times
\Bbb R$, we consider
\begin{align*}
&\Phi (q; \varphi_\ell, \zeta_\ell, \alpha_\ell)\;  :=\; 
\Psi \left (\,\omega (\ell -1) +
(\sqrt{-1}/2\pi)q^{\ell}\partial\bar{\partial}\varphi_\ell ,  \; q\, \right )\, -\,  
\\
&\; \; \qquad\qquad (C_{q,\ell -1}+ \alpha_\ell q^{\ell
+1})\, \biggl \{ 1 -\left (\,\mathcal{Y}(\ell -1)/q + q^{\ell + 1}\mathcal{Y}_{\ell}\, \right )
\cdot \log \{h(\ell -1)\exp (-q^{\ell}\varphi_{\ell} ) \} \\
&\quad\qquad\qquad  + \;\; R \left (\mathcal{Y}(\ell -1)/q + q^{\ell + 1}\mathcal{Y}_{\ell}; \;
h(\ell -1)\exp (-q^{\ell}\varphi_{\ell} ) \right ) \biggl \}.
\end{align*}
By the induction hypothesis,  $\Xi (\omega (\ell -1), \mathcal{Y}(\ell - 1), C_{q, \ell -1})
\equiv 0$ modulo $q^{\ell + 1}$. Since $\Phi (q; 0,0,0) = 
\Xi (\omega (\ell -1), \mathcal{Y}(\ell -
1), C_{q, \ell -1})$, we have 
$$
\Phi (q; 0,0,0) \; \equiv \; 
u_{\ell} q^{\ell +1},
\qquad \text{modulo $q^{\ell +2}$,}
$$
for some real-valued 
$K$-invariant smooth function $u_{\ell}$ on $M$.
Let $(\varphi_\ell,  \zeta_\ell, \alpha_k)
\in \operatorname{Ker} \mathcal{D}_0^{\perp}\times \operatorname{Ker} \mathcal{D}_0
\times\Bbb R$. Since $\varphi_k$ is $K$-invariant, by $\mathcal{V}\in \frak k$, we see that
$\sqrt{-1}\,\mathcal{V}\,\varphi_k$ is a real-valued function on $M$.
Note also that $\mathcal{Y}(0) = (\sqrt{-1}\mathcal{V}/2)\,q^2$.
Then  the variation formula for the scalar curvature (see for instance
\cite{Ca}, (2.5)) shows that, modulo $q^{\ell +2}$,
\begin{align*}
 & \Phi (q; \varphi_\ell, \zeta_\ell, \alpha_\ell) \\
&\equiv 
\Phi (q; 0,0,0 )
+ \frac{q^{\ell +1}}{2}\, (- \mathcal{D}_0 + \sqrt{-1}\,\mathcal{V} )\varphi_{\ell}
-\alpha_{\ell}q^{\ell +1}+ q^{\ell +1} h_0^{-1}(\mathcal{Y}_{\ell}\cdot h_0)
- \frac{\sqrt{-1}}{2}\,\mathcal{V}\,\varphi_{\ell}\,q^{\ell +1}\\
&\equiv  \left \{ u_{\ell} - \mathcal{D}_{0}(\varphi_{\ell}/2) - \alpha_{\ell} - 
\hat{F}_m(\mathcal{Y}_{\ell}) + e_0^{-1}(\sqrt{-1}\,\mathcal{Y}_{\ell})
\right \}\,q^{\ell + 1},
\end{align*}
where we put $\hat{F}(\mathcal{Y} ) := \{ c_1(L)^n [M]\}^{-1} 
2\pi  F(\sqrt{-1}\,\mathcal{Y})$ for each $\mathcal{Y} \in \sqrt{-1}\, \frak z$.
By setting $\mu_{\ell}:= \{c_1(L)^n[M]\}_{}^{-1}(\int_M u_{\ell}\omega_0^n)$,
we write $u_{\ell}$ as a sum 
$$
{u}_{\ell} \;=\; \mu_{\ell} + u'_{\ell}  + u''_{\ell},
$$
 where
$u'_{\ell} := (1- P)({u}_{\ell}- \mu_{\ell})\in \operatorname{Ker} \mathcal{D}_0^{\perp}$ and
$u''_{\ell} :=  P({u}_{\ell}- \mu_{\ell})\in \operatorname{Ker} \mathcal{D}_0$.
Now, let $\varphi_{\ell}$ be the unique element of $\operatorname{Ker}
\mathcal{D}_0^{\perp}$ such that  $\mathcal{D}_0(\varphi_{\ell} /2) = u'_{\ell}$. Moreover, we put
$$
\zeta_{\ell} \,:=\, u''_{\ell}\quad \text{ and }\quad
\alpha_{\ell}\, :=\, \mu_{\ell} - \hat{F}(\mathcal{Y}_{\ell}).
$$
Then by $\mathcal{Y}_{\ell} = \sqrt{-1} e_0(\zeta_{\ell} )
= \sqrt{-1}\,e_0(u''_{\ell})$, we obtain
\begin{align*}
 &\Phi (q; \varphi_\ell, \zeta_\ell, \alpha_\ell) \;\equiv \;
\left \{ \mu_{\ell} + u'_{\ell} + u''_{\ell} - \mathcal{D}_{0}(\varphi_{\ell}/2) - \alpha_{\ell} - 
\hat{F}_m(\mathcal{Y}_{\ell}) + e_0^{-1}(\sqrt{-1}\,\mathcal{Y}_{\ell})
\right \}\,q^{\ell + 1}\\
&\equiv \; \{\,u''_{\ell} + e_0^{-1}(\sqrt{-1}\,\mathcal{Y}_{\ell})\,\}\,q^{\ell + 1}\;
\equiv\; 0, \qquad \text{mod $q^{\ell +2}$,}
\end{align*}
as required.
Write $\sqrt{-1}\,\mathcal{V}/2$ as $\mathcal{Y}_0$ for simplicity.
Now, for the real Lie subalgebra $\frak b$ of $\frak z$
generated by $\mathcal{Y}_k$, $k =0,1,2,\dots$,   its
complexification $\frak b^{\Bbb C}$
 in $\frak z^{\Bbb C}$ generates a complex Lie subgroup
$B^{\Bbb C}$ of $Z^{\Bbb C}$.
Then it is easy to check that the algebraic subtorus
$T$ of $Z^{\Bbb C}$ obtained as
the closure of $B^{\Bbb C}$ in  $Z^{\Bbb C}$
has the required properties.

\medskip
{\em Remark $6.5$}. In Theorem C, assume that $\omega_0$ is a K\"ahler metric of
constant scalar curvature, and moreover that the actions
$\rho_{m(\nu )}$, $\nu = 1,2,\dots,$ coincide
for all sufficiently large $\nu$. Then by
\cite{M4}, the trivial group  $\{1\}$ can be chosen as the algebraic subtorus $T$ 
above of $Z^{\Bbb C}$.

\end{document}